\documentclass[12pt]{article}

\usepackage{amsfonts}
\usepackage{amsmath}
\usepackage[all]{xy}
\newtheorem{defn}{Definition}

\newtheorem{thm}{Theorem}

\newtheorem{cor}{Corollary}

\newtheorem{ex}{Example}
\newtheorem{rem}{Remark}

\title{A partition of connected graphs}

\author{Gus Wiseman\\
\small Department of Mathematics,
\small University of California,
\small One Shields Ave.,
\small Davis, CA 95616\\
\small \texttt{gus@nafindix.com}
}

\date{\small 
Published in The Electronic Journal of Combinatorics, Jan 7, 2005\\
\small Mathematics Subject Classifications: 05C30, 05C05}

\begin{document}
\bibliographystyle{plain}
\maketitle

\begin{abstract}
We define an algorithm $k$ which takes a connected graph $G$ on a totally ordered vertex set and returns an increasing tree $R$ (which is not necessarily a subtree of $G$). We characterize the set of graphs $G$ such that $k(G)=R$. Because this set has a simple structure (it is isomorphic to a product of non-empty power sets), it is easy to evaluate certain graph invariants in terms of increasing trees. In particular, we prove that, up to sign, the coefficient of $x^q$ in the chromatic polynomial $\chi_G(x)$ is the number of increasing forests with $q$ components that satisfy a condition that we call $G$-connectedness. We also find a bijection between increasing $G$-connected trees and broken circuit free subtrees of $G$.
\end{abstract}

We will work with finite labeled simple graphs. Usually we will identify a graph $G$ with its edge set; this should not cause any serious ambiguities. If the vertex set is $V$ then we say that $G$ is a graph on $V$. A (spanning) subgraph $Q$ of $G$ is a graph with the same vertex set as $G$ and a subset of the edges of $G$. The notation $Q\subseteq G$ means $Q$ is a subgraph of $G$. A rooted graph is a graph with a distinguished vertex called the root.

Define $\mathrm{link}(v,S)$ to be the set of all possible edges joining $v$ to an element of $S$ (so if $v\notin S$, $\mathrm{link}(v,S)$ has $|S|$ elements). If $G$ is a graph on $V$ and $S\subseteq V$, we define the restriction of $G$ to $S$, $G|_S$, to be the graph on $S$ whose edge set consists of all edges of $G$ with both ends in $S$.

We will use the symbols $\pi$ and $\sigma$ to denote set partitions. The notation $\pi\vdash S$ means $\pi$ is a set partition of the set $S$. The length (number of blocks) of $\pi$ is denoted by $\ell(\pi)$. A set partition $\sigma$ is called a refinement of a set partition $\pi$ if every block of $\sigma$ is contained in some block of $\pi$.

To each graph $G$ on $V$ there corresponds a set partition $s(G)$ such that two vertices $v,w\in V$ are in the same block of $s(G)$ if and only if there is a path in $G$ from $v$ to $w$. Equivalently, $s(G)$ is the maximal set partition of $V$ whose blocks are connected. The restriction of $G$ to a block of $s(G)$ is called a component of $G$.

If $G$ is a rooted connected graph on $V$ with root $r$, we will call the set partition $\pi=s(G|_{V-\{r\}})$ of $V-\{r\}$ the depth-first partition of $G$. To obtain a connected subgraph of a rooted connected graph $G$ on $V$, we can choose, for each block $\pi_i$ of $\pi$, a connected subgraph of $G|_{\pi_i}$ and a nonempty set of edges (in $G$) connecting $r$ to $\pi_i$. In fact, every connected subgraph of $G$ can be obtained in this way. Our Theorem~\ref{thm:nm} may be regarded as an iteration of this correspondence. The depth-first partition and this correspondence have been studied by Gessel~\cite{MR96f:05092}. 

A forest is a graph with no circuits. A tree is a connected forest. A basic property of trees is that there is a unique path (a sequence of distinct, adjacent vertices) between any two vertices. The distance between two vertices is defined to be the length of this path. In a rooted tree, the height of a vertex is defined to be its distance from the root. A vertex $w$ is called a descendant of a vertex $v$ (or $v$ is called an ancestor of $w$) if the heights of the vertices on the unique path from $v$ to $w$ are increasing (so in particular $v$ is always a descendant of itself). We define the join of $v$ and $w$ to be their unique common ancestor on the unique path between them.

Let $R$ be a rooted tree on the vertex set $V$, and let $v\in V$. We define $\mathrm{des}(v,R)\subseteq V$ to be the set of descendants of $v$ (including $v$). If $v$ is not the root of $R$, we define $\mathrm{parent}(v,R)\in V$ to be the closest vertex to $v$ in $R$ which is not a descendant of $v$. A rooted tree is increasing (according to a total order on $V$) if for each $v\in V$ and $w\in\mathrm{des}(v,R)$ we have $v\leq w$. Consequently, the root of an increasing tree must be the smallest element of $V$.

\begin{defn}
Let $R$ be a rooted tree on the totally ordered vertex set $V$ with root $r$, and let $v\in V-\{r\}$. Define $J(v,R)=\mathrm{link}(\mathrm{parent}(v,R),\mathrm{des}(v,R))$. If $G$ is a graph on $V$ and if for each $v\in V-\{r\}$ we have $J(v,R)\cap G\neq\emptyset$ then we say that $R$ is $G$-connected.
\end{defn}

Note that the sets $J(v,R)$ (as $v$ ranges over $V-\{r\}$) are disjoint. Also note that a $G$-connected tree need not be a subgraph of $G$ and that $G$ must be connected for any rooted tree to be $G$-connected.

\begin{defn}\label{defn:k}
For each connected graph $G$ on a totally ordered vertex set $V$, define an increasing $G$-connected tree $k(G)$ by the following algorithm:
\begin{enumerate}
\item Let $H$ be an empty graph on $V$, and set $S=V$.
\item Let $\pi$ be the depth-first partition of $G|_S$ rooted at r=the smallest vertex in $S$. Add edges to $H$ connecting $r$ to the smallest vertex in each block of $\pi$.
\item For each block $\pi_i$ of $\pi$ with more than one element, return to step 2 with $S=\pi_i$.
\item Return $k(G)=H$.
\end{enumerate}
\end{defn}

\begin{ex}
The $6$ increasing trees on $V=\{1,2,3,4\}$ are listed vertically. To the right of each increasing tree $R$ are listed the subtrees $T$ of the complete graph on $V$ such that $k(T)=R$ (we have omitted the $22$ connected subgraphs which are not trees). The breaks are indicated by dotted lines (see Theorem~\ref{thm:bc}).
\\

$
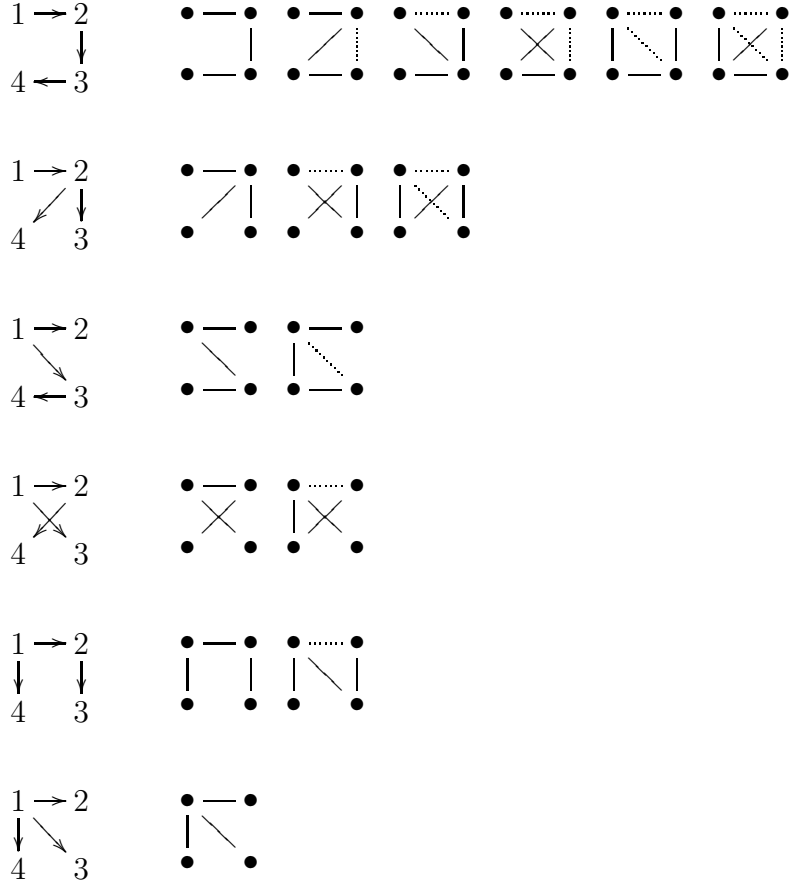
\xymatrix @-1pc{
1 \ar@{->}[r] & 2 \ar@{->}[d] \\
4 & 3 \ar@{->}[l] 
}
$
\qquad
$
\xymatrix @-1pc{
\bullet \ar@{-}[r] & \bullet \ar@{-}[d] \\
\bullet & \bullet \ar@{-}[l] 
}
$
$
\xymatrix @-1pc{
\bullet \ar@{-}[r] & \bullet \ar@{-}[dl] \ar@{.}[d]\\
\bullet \ar@{-}[r] & \bullet 
}
$
$
\xymatrix @-1pc{
\bullet \ar@{-}[dr] \ar@{.}[r] & \bullet \ar@{-}[d] \\
\bullet & \bullet \ar@{-}[l] 
}
$
$
\xymatrix @-1pc{
\bullet \ar@{-}[dr] \ar@{.}[r] & \bullet \ar@{-}[dl] \ar@{.}[d]\\
\bullet & \bullet \ar@{-}[l] 
}
$
$
\xymatrix @-1pc{
\bullet \ar@{-}[d] \ar@{.}[r] \ar@{.}[dr]& \bullet \ar@{-}[d] \\
\bullet & \bullet \ar@{-}[l] 
}
$
$
\xymatrix @-1pc{
\bullet \ar@{-}[d] \ar@{.}[r] \ar@{.}[dr]& \bullet \ar@{-}[dl] \ar@{.}[d]\\
\bullet & \bullet \ar@{-}[l] 
}
$
\\ \\

$
\xymatrix @-1pc{
1 \ar@{->}[r] & 2 \ar@{->}[d] \ar@{->}[dl] \\
4 & 3
}
$
\qquad
$
\xymatrix @-1pc{
\bullet \ar@{-}[r] & \bullet \ar@{-}[d] \\
\bullet \ar@{-}[ur] & \bullet  
}
$
$
\xymatrix @-1pc{
\bullet \ar@{-}[dr] \ar@{.}[r] & \bullet \ar@{-}[ld] \\
\bullet & \bullet \ar@{-}[u]
}
$
$
\xymatrix @-1pc{
\bullet \ar@{-}[d] \ar@{.}[r] \ar@{.}[dr]& \bullet \ar@{-}[d] \\
\bullet \ar@{-}[ur] & \bullet  
}
$
\\ \\

$
\xymatrix @-1pc{
1 \ar@{->}[r] \ar@{->}[dr] & 2\\
4 & 3 \ar@{->}[l]
}
$
\qquad
$
\xymatrix @-1pc{
\bullet \ar@{-}[r] \ar@{-}[dr] & \bullet \\
\bullet \ar@{-}[r] & \bullet  
}
$
$
\xymatrix @-1pc{
\bullet \ar@{-}[r] \ar@{-}[d] \ar@{.}[dr] & \bullet \\
\bullet & \bullet \ar@{-}[l]
}
$
\\ \\

$
\xymatrix @-1pc{
1 \ar@{->}[r] \ar@{->}[dr] & 2 \ar@{->}[dl] \\
4 & 3
}
$
\qquad
$
\xymatrix @-1pc{
\bullet \ar@{-}[r] \ar@{-}[dr] & \bullet \\
\bullet \ar@{-}[ur] & \bullet  
}
$
$
\xymatrix @-1pc{
\bullet \ar@{-}[dr] \ar@{-}[d] \ar@{.}[r] & \bullet \\
\bullet \ar@{-}[ur] & \bullet
}
$
\\ \\

$
\xymatrix @-1pc{
1 \ar@{->}[d] \ar@{->}[r] & 2 \ar@{->}[d]\\
4 & 3
}
$
\qquad
$
\xymatrix @-1pc{
\bullet \ar@{-}[d] \ar@{-}[r] & \bullet \\
\bullet & \bullet \ar@{-}[u]
}
$
$
\xymatrix @-1pc{
\bullet \ar@{-}[dr] \ar@{-}[d] \ar@{.}[r] & \bullet \\
\bullet & \bullet \ar@{-}[u]
}
$
\\ \\

$
\xymatrix @-1pc{
1 \ar@{->}[r] \ar@{->}[dr] \ar@{->}[d]& 2\\
4 & 3
}
$
\qquad
$
\xymatrix @-1pc{
\bullet \ar@{-}[d] \ar@{-}[r] & \bullet \\
\bullet & \bullet \ar@{-}[ul]
}
$
\end{ex}

There is a different algorithm, called depth-first search, which produces subforests of $G$. Some enumerative applications of this algorithm have been studied by Gessel and Sagan~\cite{MR97d:05149}. A distinguishing difference between depth-first search and our algorithm is that depth-first search only follows the edges of $G$, whereas here we add edges connecting to the smallest vertex in each block of $\pi$ regardless of whether these are edges of $G$. The algorithms are related in that if $G$ is a connected graph and $R$ is a depth-first search subtree of $G$ then parts 2 and 3 of the next theorem hold (although the converse is not true).

\begin{thm}\label{thm:nm}Let $G$ be a connected graph on a totally ordered vertex set $V$, and let $R$ be an increasing $G$-connected tree on $V$. Then the following are equivalent:
\begin{enumerate}
\item $k(G)=R$
\item For each vertex $v\in V$, $G|_{\mathrm{des}(v,R)}$ rooted at $v$ is connected and has the same depth-first partition as $R|_{\mathrm{des}(v,R)}$ rooted at $v$.
\item For each non-root vertex $v\in V-\{r\}$ there is a nonempty set $E(v)\subseteq J(v,R)$ such that $G=\bigcup_{v\in V-\{r\}}E(v)$.
\end{enumerate}
\end{thm}

\emph{Proof.} $1\Leftrightarrow2$ This follows easily from Definition~\ref{defn:k}.

$2\Rightarrow3$ Let $E(v)=J(v,R)\cap G$. We need to show that every edge of $G$ lies in some $E(v)$. Let $e\in G$ and let $v<w$ be the vertices of $e$. We will show that $w$ is a descendant of $v$. Suppose this is false, and let $u$ be their join. Then $e\in G|_{\mathrm{des}(u,R)}$, so $v$ and $w$ are in the same block of the depth first partition of $G|_{\mathrm{des}(u,R)}$. This is a contradiction because they are in different blocks of the depth first partition of $R|_{\mathrm{des}(u,R)}$. Now, since $w$ is a descendant of $v$, there is a unique vertex $z\in V$ (possibly equal to $w$) such that $\mathrm{parent}(z)=v$ and $w\in\mathrm{des}(z)$. Hence $e\in J(z,R)\cap G$.

$3\Rightarrow2$ This is certainly true if $v$ (in part 2) is a leaf of $R$ (its only descendant is itself). Let $v\in V$ and suppose it is true for all $w\in\mathrm{des}(v,R)-\{v\}$. Let $\pi$ be the depth-first partition of $R|_{\mathrm{des}(v,R)}$. Then $G|_{\pi_i}$ is connected by the inductive hypothesis. Furthermore, $G$ contains an edge connecting $v$ to $\pi_i$ because $\pi_i$ contains a vertex $w$ whose parent in $R$ is $v$ and $J_G(w,R)$ consists of edges connecting $v$ to $\pi_i$. Hence $G|_{\mathrm{des}(v,R)}$ is connected. Clearly $\pi$ is a refinement of the depth-first partition of $G|_{\mathrm{des}(v,R)}$ (because $G|_{\pi_i}$ is connected), so to show that they are equal we have only to show that if $x$ and $y$ are in different blocks of $\pi$ then they are in different blocks of the depth-first partition of $G$. Let $x<y\in V$ be in different blocks of $\pi$, and suppose $G$ has an edge between $x$ and $y$. Then $y$ is a descendant of $x$ in $R$ because every edge of $J_G(w,R)$ (for any $w\in V$) connects a vertex to one of its descendants. This contradicts the fact that they are in different blocks of the depth-first partition of $R|_{\mathrm{des}(v,R)}$. $\square$

\begin{rem}
Actually the condition in Theorem~\ref{thm:nm} that $R$ be $G$-connected is not necessary because if $R$ is not $G$-connected then parts 1, 2 and 3 will be false.
\end{rem}

Some algebraic invariants of graphs can be simply expressed in terms of connected subgraphs. We can use the algorithm $k$ to express such invariants in terms increasing trees. Moreover, Theorem~\ref{thm:nm} shows that the set $k^{-1}(R)$ has a simple structure, as illustrated by the next theorem.

\begin{defn}Let $G$ be a connected graph on $V$. Define
\[
\eta^G(t)=\sum_{\substack{Q\subseteq G \\ \mathrm{connected}}}t^{|Q|}
\]
where $|Q|$ denotes the number of edges in $Q$.
\end{defn}

\begin{thm}\label{thm:polyv}
\[
\eta^G(t)=\sum_{\substack{R \\ \mathrm{increasing} \\ G-\mathrm{connected}}}
\prod_{v\in V-\{r\}}
[(1+t)^{|J(v,R)\cap G|}-1]
\]
\end{thm}

\emph{Proof.} We have
\[
\eta^G(t)=\sum_{\substack{R \\ \mathrm{increasing} \\ G-\mathrm{connected}}}\sum_{\substack{Q\subseteq G \\ k(Q)=R}}t^{|Q|}
\]
Now, the generating function for the cardinality of nonempty subsets of a set $S$ is
\[
f_S(x)=\sum_{\emptyset\neq T\subseteq S}x^{|T|}=(1+x)^{|S|}-1
\]
Hence from Theorem~\ref{thm:nm} part 3,
\[
\sum_{\substack{Q\subseteq G \\ k(Q)=R}}t^{|Q|}=\sum_{\substack{Q=\bigcup_{v\in V-\{r\}}E(v) \\ \emptyset\neq E(v)\subseteq J(v,R)\cap G}}t^{|Q|}=
\prod_{v\in V-\{r\}}f_{J(v,R)\cap G}(t)
\]
from which the result follows. $\square$

The chromatic polynomial $\chi_G(x)$ of a graph $G$ is a polynomial which evaluates to the number of proper colorings of $G$ with $x$ colors. The subgraph expansion of $\chi_G(x)$ is
\[
\chi_G(x)=\sum_{Q\subseteq G}(-1)^{|Q|}x^{c(Q)}
\]
where $c(Q)$ is the number of components of $Q$. See~\cite{MR95h:05105} for background on the chromatic polynomial.

We define an increasing $G$-connected forest $R$ to be a forest where each component $R|_{s(R)_i}$ is an increasing $G|_{s(R)_i}$-connected tree. For a graph $G$, let $t(G)$ be the (integer) partition whose parts are the sizes of the blocks of $s(G)$. For background on the chromatic symmetric function $X_G=X_G(x_1,x_2,\ldots)$ of a graph $G$, see~\cite{MR96b:05174} and~\cite{MR2000c:05152}. For background on the chromatic symmetric function in non-commuting variables $Y_G=Y_G(x_1,x_2,\ldots)$, see~\cite{MR2002d:05124}.

\begin{cor}\label{cor:inv}Let $G$ be a graph on a totally ordered vertex set $V$ with $|V|=n$.
\begin{enumerate}
\item The coefficient of $(-1)^{n-1}x$ in the chromatic polynomial $\chi_G(x)$ is the number of increasing $G$-connected trees.
\item The coefficient of $(-1)^{n-q}x^q$ in the chromatic polynomial $\chi_G(x)$ is the number of increasing $G$-connected forests with $q$ components (or, equivalently, with $n-q$ edges).
\item The coefficient of $(-1)^{n-\ell(\lambda)}p_\lambda$ in the chromatic symmetric function $X_G$ is the number of increasing $G$-connected forests $R$ such that $t(R)=\lambda$.
\item The coefficient of $(-1)^{n-\ell(\pi)}p_\pi$ in the chromatic symmetric function in non-commuting variables $Y_G$ is the number of increasing $G$-connected forests $R$ such that $s(R)=\pi$.
\end{enumerate}
\end{cor}

\emph{Proof. } 1. Let $a^G$ be the coefficient of $x$ in $\chi_G(x)$. From the subgraph expansion we have
\[
a^G=\sum_{\substack{Q\subseteq G \\ \mathrm{connected}}}(-1)^{|Q|}=\eta^G(-1)=\sum_{\substack{R \\ \mathrm{increasing} \\ G-\mathrm{connected}}}\prod_{v\in V-\{r\}}(-1)
\]

We don't need to worry about $0^0$ because the $G$-connectedness of $R$ implies that $J(v,R)\cap G$ is never empty.

4. We will prove part 4, the others being simple specializations. Let $H_\pi^G$ be the number of increasing $G$-connected forests $R$ such that $s(R)=\pi$, and let $H^G$ be the number of increasing $G$-connected trees. Then using part 1 we have
\begin{equation}\label{eq:an}
H_\pi^G=\prod_{i=1}^{\ell(\pi)}H^{G|_{\pi_i}}=(-1)^{n-\ell(\pi)}\prod_{i=1}^{\ell(\pi)}\sum_{\substack{Q\subseteq G|_{\pi_i} \\ connected}}(-1)^{|Q|}
\end{equation}
The subgraph expansion of $Y_G$ is
\[
Y_G=\sum_{Q\subseteq G}(-1)^{|Q|}p_{s(Q)}
\]
Hence
\[
Y_G=\sum_{\pi\vdash V} p_\pi\sum_{\substack{Q\subseteq G \\ s(Q)=\pi}}(-1)^{|Q|}=
\sum_{\pi\vdash V} p_\pi\prod_{i=1}^{\ell(\pi)}\sum_{\substack{Q\subseteq G|_{\pi_i} \\ connected}}(-1)^{|Q|}
\]
Substituting~\eqref{eq:an}, we obtain the desired result. $\square$

If $G$ is a graph on a totally ordered vertex set $V$, we extend the ordering of the vertices to an ordering of the edges lexicographically. A broken circuit of $H\subseteq G$ is a set of edges $B\subseteq H$ such that there is some edge $e\in G$, smaller than every edge of $B$, such that $B\cup e$ is a circuit. Note that $B$ being a broken circuit of $H$ depends both on $H$ and $G$. If $H\subseteq G$ contains no broken circuits then it is called broken circuit free. Note that if $H$ contains a circuit then it also contains a broken circuit. Consequently, a broken circuit free subgraph is always a forest. If $T\subseteq G$ is a subtree of $G$ and the edge $e\in G$, $e\notin T$ is the smallest edge in the unique circuit in $T\cup\{e\}$ then we will call $e$ a break in $T$. Hence the set of breaks in a subtree $T$ is in bijection with the set of broken circuits of $T$.

Whitney's Broken Circuit Theorem~\cite{logicalexpansion} shows that if $G$ is a connected graph with $n$ vertices, the coefficient of $(-1)^{n-1}x$ in $\chi_G(x)$ is the number of broken circuit free subtrees of $G$. Hence there should be a bijection between broken circuit free subtrees and increasing $G$-connected trees. 

\begin{thm}\label{thm:bc}Let $V$ be a totally ordered vertex set with smallest element $r$, and let $G$ be a connected graph on $V$. Let $T\subseteq G$ be a subtree of $G$, and let $R=k(T)$. Let $E(v)$ for $v\in V-\{r\}$ be as in Theorem~\ref{thm:nm} part 3. Then $E(v)$ contains only one element $e(v)$ (otherwise $T$ would have more than $|V|-1$ edges so it could not be a tree). For $v\in V-\{r\}$, let $d(v)$ be the set of elements of $J(v,R)\cap G$ which are smaller than $e(v)$. Then the set of breaks in $T$ is
\[
\bigcup_{v\in V-\{r\}}d(v)
\]
\end{thm}

\emph{Proof.} Let $J=\bigcup_{v\in V-\{r\}}J(v,R)\cap G$. Since $k(G)$ may be different from $R$, $J$ may be different from $G$. We will first show that if $e\in G$ but $e\notin J$ then $e$ is not a break. Let $v<w\in V$ be the vertices of $e$. Then $w$ is not a descendant of $v$ because otherwise we would have $e\in J$. Let $u\in V$ be the join of $v$ and $w$ in $R$. Then Theorem~\ref{thm:nm} part 2 implies that $u$ is also the join of $v$ and $w$ in $T|_{\mathrm{des}(u,R)}$ (rooted at $u$). Therefore, the cycle created by adding $e$ to $T$ contains an edge connected to $u$. Since $u<v<w$, $e$ cannot be a break.

Now suppose $e\in J(v,R)\cap G$ is smaller than $e(v)$. We will show that $e$ is a break. Let $H=T|_{\mathrm{des}(v,R)\cup\mathrm{parent}(v,R)}$. Then $\mathrm{parent}(v,R)$ is the smallest vertex in the vertex set of $H$. Therefore, $e$ is smaller than any other edge in $H$. Since $H$ is a tree, adding $e$ would create a unique circuit in $H$. Hence $e$ is a break.

Now suppose $e\in J(v,R)\cap G$ is larger than $e(v)$. Then, letting $H$ be as before, we see that $e(v)$ must belong to the circuit which $e$ creates. But $e(v)$ is smaller than $e$, so $e$ cannot be a break. $\square$

\begin{cor}\label{cor:three}The function
\[
f(R)=\bigcup_{v\in V-\{r\}}\mathrm{min}(J(v,R)\cap G)
\]
is a bijection between increasing $G$-connected trees and broken circuit free subtrees, and $f^{-1}(T)=k(T)$.
\end{cor}

Of course, this bijection generalizes to a bijection between increasing $G$-connected forests with $q$ components and broken circuit free subforests of $G$ with $q$ components.

\bibliography{refs}

\end{document}